

\magnification=\magstep1
\baselineskip=16 pt
\mathsurround=2 pt


\def\gr{\it}
\def\tit#1{\bigbreak\noindent{\bf #1}\medskip\nobreak}
\def\PF{\noindent{\it Proof:\/} }
\def\FP{\hbox{\quad\vrule height6pt width4pt depth0pt}}

\def\IFF{if and only if}
\def\AND{{\rm\ and\ }}
\def\OR{{\rm\ or\ }}
\def\IMP{\Rightarrow}

\def\0{\hbox{{\rm\O}}}
\def\sc{\subset}

\def\S#1{\Sigma(#1)}
\def\F{{\cal F}}
\def\<#1>{\langle#1\rangle}
\def\EXT#1{\,\overline{#1}\,}
\def\K#1{{\rm K}(#1)}

\def\R{{\bf R}}
\def\N{{\bf N}}

\def\A{{\rm A}}
\def\D{{\rm D}}
\def\U{{\rm U}}

\def\SR{\S\R}
\def\FA{\F(\A)}
\def\FR{\F(\R)}
\def\KFR{{\rm K}[\FR]}

\def\id#1{{\rm id}_{#1}}

\def\:{\colon}
\def\o{\circ}
\def\f{f}
\def\g{g}

\def\RR{\EXT\R}
\def\NN{\EXT\N}

\def\UU{\EXT\U}
\def\FFR{\EXT{\FR}}
\def\FRR{\F(\RR)}
\def\SFR{\Sigma[\FR]}

\def\ff{\EXT\f}
\def\xx{\EXT x}

\def\a{\alpha}
\def\b{\beta}
\def\ga{\gamma}
\def\k{\kappa}
\def\m{\mu}
\def\n{\nu}
\def\x{\xi}
\def\t{\tau}
\def\y{\upsilon}

\def\8{\infty}
\def\d{\partial}

\def\dom{{\rm dom}}
\def\im{{\rm im}}
\def\={\approx}
\def\@{\sim}


\hrule height0pt
\vfil
{\bf VIRTUAL CALCULUS --- PART II}
\bigskip\bigskip\bigskip
{\bf S\'ergio F.~Cortizo}
\bigskip\bigskip
Instituto de Matem\'atica e Estat\'\i stica, Universidade de S\~ao
Paulo

Cidade Universit\'aria, Rua do Mat\~ao, 1010

05508-900, S\~ao Paulo, SP, Brasil
\smallskip
cortizo@ime.usp.br
\bigskip\bigskip\bigskip
{\bf Abstract}
\bigskip
\item{}A simultaneous extension of real numbers set and the class of real
functions is discussed.
\smallskip
PACS \ 02.90.+p
\vfil\eject

\tit{I. Introduction}

In Ref.~1 we presented a process of extending which can be applied to any
set, including several sets simultaneously. In the first part of this
work$^2$, we applied this process to the ordered field of real numbers~$\R$,
thus obtaining the set $\RR$ of {\gr virtual numbers}, which contains
infinitesimal and infinite quantities. That extension was then used in a
reorganization of Infinitesimal Calculus.

We will now apply that extending process simultaneously to~$\R$ and to the
set of real functions, thus defining {\gr virtual functions\/} as classes of
sequences of real functions. By means of some identifications, we can
consider these mathematical objects as maps between virtual numbers, i.e., as
functions from subsets of~$\RR$ into~$\RR$.

Having done this, we can extend the constructions and techniques of
Infinitesimal Calculus to the virtual functions. For example, we can define
its derivatives and so prove that they can be calculated according to the
usual derivation rules. We can also define the integration of a virtual
function between two virtual numbers, which provides another virtual number
as a result, and then generalize the Fundamental Theorem of Calculus to those
integrals.

This procedure ends up broadening the reach of Calculus techniques,
systematizing many ideas which have been used for the last decades in
specific contexts. We can mention, for example, Delta Calculus created by
Dirac$^3$, which has been object of various attempts of mathematical rigorous
formalization. This application of virtual functions will be presented in a
subsequent work.

In Sec.~II we discuss the simultaneous virtual extension of real numbers and
functions. In Sec.~III, we show that the virtual objects thus obtained can be
manipulated practically as if they were real. In Secs. IV, V and~VI we extend
the basic constructions of Calculus to the virtual functions. Sec.~VII is
dedicated to the Fundamental Theorem of Calculus. In Sec.~VIII we introduce
the concept of virtual sequence, and point out other possible applications of
the virtual extension process to Calculus.

\tit{II. Virtual Functions}

For any set~$\A$, we will represent the class of all functions whose domain
and image are subsets of~$\A$ by~$\FA$. In other words, $\FA$ is the set of
functions $\f\:\D\to\A$ such that $\D\subset\A$. We will consider the {\gr
empty relation\/}~$\0$ as a member of~$\FA$, for each set~$\A$ (the domain
and the image of this function are both empty).

This way, $\FR$ denotes the set of {\gr real functions}, i.e., the set of all
functions whose domain and image are subsets of the ordered field $\R$ of
real numbers. The empty function $\0\in\FR$ should not be mistaken by the
null real function, which constantly equals zero and whose domain is the
whole real line.\smallskip

We will now consider the simultaneous virtual extension of~$\R$ and~$\FR$,
i.e., we will apply the process of virtual extension to the disjoint union
$\U$ of these two sets$^1$. The members of the virtual extension $\RR$ of the
subset $\R\subset\U$ will be called {\gr virtual numbers\/}, and the members
of the virtual extension $\FFR$ of the subset $\FR\subset\U$ will be called
{\gr virtual functions}.\smallskip

A generic element of the extension $\UU$ is an equivalence class of sequences
on~$\U$, i.e., of sequences formed either by real numbers or real functions.
The members of~$\RR$ are the classes of sequences which end taking only
values in~$\R$, whereas the virtual functions are the classes of sequences
which end taking only values in~$\FR$.\smallskip

We will consider $\R\sc\RR$ according to the identification $\R=\K\R$, and
$\FR\sc\FFR$ according to the identification $\FR=\KFR$, i.e., the element
$x\in\R$ will be identified with the class $\<x,x,\ldots>$ of the constant
sequence at~$x$, and the function $\f\in\FR$ will be identified with the
class $\<\f,\f,\ldots>$ of the constant sequence at~$\f$. Thus, we will omit
both the ``bar'' which distinguishes $x\in\R$ from~$\xx\in\K\R$ and that
which distinguishes $\f\in\FR$ from~$\ff\in\KFR$. In this manner, any real
number can be understood as a virtual number, and every real function can be
understood as a virtual function. Nevertheless, it is clear that the
reciprocal statements do not hold.\smallskip

Throughout this whole work, real numbers and functions will be denoted by
low-case Latin letters ($x$, $a$, $\f$, $\g,\ldots$), whereas virtual numbers
and functions will be denoted by low-case Greek letters ($\x$, $\a$, $\phi$,
$\psi,\ldots$). The letter `$\pi$' is an exception keeping its usual
mathematical meaning: $\pi\in\R$ is the constant ratio between the
circumference and the diameter of a circle.\smallskip

A virtual function $\phi\in\FFR$, as defined above, is not a relation between
elements of two sets. Therefore, it has neither a {\gr domain\/} nor an {\gr
image\/} according to the usual meaning of these terms. In spite of that, we
can consider the virtual extension of the relation ``is defined at'' between
real functions and real numbers, which is a relation between virtual
functions and virtual numbers. That extension allows us to decide when
$\phi\in\FFR$ {\gr is defined at\/}~$\x\in\RR$. So, we will say that the {\gr
domain\/} of a virtual function $\phi\in\FFR$ is the set of all $\x\in\RR$
which satisfy that condition. This set will be denoted
by~$\dom\phi\subset\RR$, in analogy with the notation $\dom\f\subset\R$ for
the domain of a function~$\f\in\FR$. Thus, if $\phi$ is the class of the
sequence $(\f_1,\f_2,\ldots)\in\SFR$ and $\x$ the class
of~$(x_1,x_2,\ldots)\in\SR$, then $\x\in\dom\phi$ \IFF\ there exists $n\in\N$
such that, for every $i>n$, the function $\f_i$ is defined at~$x_i$.

The domain of the virtual empty function $\0=\<\0,\0,\0,\ldots>\in\FFR$ is
the empty subset $\0\subset\RR$. It is important to note that, besides this
one, there are many other virtual functions with an empty domain. For
instance, if $\id\R\in\FR$ is the {\gr identity function\/} on~$\R$, then the
sequence of functions:
$$
    (\id\R,\0,\id\R,\0,\ldots)\in\SFR
$$
represents a virtual function, different from~$\0\in\FFR$ whose domain is
empty. It is easy to verify that the domain of
$\phi=\<\f_1,\f_2,\ldots>\in\FFR$ is empty when, for every $n\in\N$, there
exists $i>n$ such that $\f_i=\0$. Also, it is not difficult to see that if
the domain of $\phi=\<\f_1,\f_2,\ldots>$ is equal to~$\RR$, then there exists
$n\in\N$ such that $\dom\f_i=\R$ for every $i>n$.

We can consider the evaluation of a real function $\f$ at a real number $x$
in its domain $\dom\f\subset\R$ as a map which sends each pair $(\f,x)$, with
$\f\in\FR$ and~$x\in\dom\f$, to a member $\f(x)\in\R$. The virtual extension
of this map defines the {\gr evaluation of a virtual function}, which allows
us to compute the {\gr image\/} by $\phi$ of any element in its domain.
Putting it another away, if $\phi=\<\f_1,\f_2,\ldots>$ is defined
at~$\x=\<x_1,x_2,\ldots>$, then $\phi(\x)\in\RR$ is the class of sequences
which end taking the values $\f_i(x_i)$.

Following the convention introduced in the beginning of this section, we will
denote the set of all functions whose domain and image are subsets of~$\RR$
by~$\FRR$, i.e., $\FRR$ is the set of functions between members of the
virtual extension of~$\R$. This set should not be confused with the virtual
extension $\FFR$ of the set of functions between members of~$\R$.

According to that notation, the above defined ``virtual evaluation'' assigns
a function $[\phi]\in\FRR$ to each virtual function $\phi\in\FFR$. The domain
of $[\phi]$ is exactly the set $\dom\phi\subset\RR$:
$$
    [\phi] \: \dom\phi \to \RR.
$$
This assignment is such that:

{\gr If $\phi$ and~$\psi$ are two distinct virtual functions, both having a
non-empty domain, then the functions $[\phi]$ and~$[\psi]$ are also
distinct.}

\PF Let $(\f_1,\f_2,\ldots)\in\SFR$ be a representative sequence
of~$\phi$, and $(\g_1,\g_2,\ldots)\in\SFR$ a representative sequence
of~$\psi$. If $\phi\ne\psi$ then, for every $n\in\N$, there exists $i>n$
such that $\f_i\ne\g_i$. Thus, for every $n\in\N$, there exist $i>n$
and~$x_i\in\R$ such that:
$$
    x_i\in\dom\f_i\ \AND\ x_i\notin\dom\g_i;\qquad\OR
$$
$$
    x_i\notin\dom\f_i\ \AND\ x_i\in\dom\g_i;\qquad\OR
$$
$$
    x_i\in\dom\f_i\ \AND\ x_i\in\dom\g_i\ \AND\ \f_i(x_i)\ne\g_i(x_i).
$$
Therefore, at least one of the three conditions below is satisfied:

(a) For every $n\in\N$, there exist $i>n$ and~$x_i\in\R$ such that
$x_i\in\dom\f_i$ and~$x_i\notin\dom\g_i$.

(b) For every $n\in\N$, there exist $i>n$ and~$x_i\in\R$ such that
$x_i\notin\dom\f_i$ and~$x_i\in\dom\g_i$.

(c) For every $n\in\N$, there exist $i>n$ and~$x_i\in\R$ such that
$x_i\in\dom\f_i$, $x_i\in\dom\g_i$ and~$\f_i(x_i)\ne\g_i(x_i)$.

If (a) holds then there exists $\x\in\dom\phi\ne\0$ which does not belong
to~$\dom\psi$. On the other hand, if (b) holds then there exists
$\x\in\dom\psi\ne\0$ which does not belong to~$\dom\phi$. Finally,
condition (c) implies the existence of a virtual number $\x$ which
belongs to the domain of both functions, and for which
$[\phi](\x)\ne[\psi](\x)$.

In any of these three cases we have $[\phi]\ne[\psi]$.\FP

Having that result in mind, we will go on by identifying all virtual
functions with an empty domain with the empty virtual function~$\0\in\FFR$.
So, we will be allowed to say:

{\gr Two virtual functions $\phi$ and~$\psi$ are equal \IFF\ they have
the same domain $\D\subset\RR$ and $[\phi](\x)=[\psi](\x)$, for each
$\x\in\D$.}

Therefore, we can see there is no significant distinction between virtual
functions $\phi\in\FFR$ and the corresponding $[\phi]\in\FRR$ anymore. So,
they will be identified to each other, not to overload notation needlessly.
In other words, we will consider:
$$
    \FFR \subset \FRR.
$$
In this sense, we will simply write:
$$\eqalign{
    \phi \:\ & \dom\phi \to \RR  \cr
             & \x \mapsto \phi(\x),\cr}
$$
and treat virtual functions as relations between virtual numbers. For
example, we define the {\gr image\/} of a virtual function $\phi$ as the set
of virtual numbers which are the image of some member of its domain. That
image will be denoted by $\im\phi\subset\RR$, as an analogy to the notation
$\im\f\subset\R$ for the image of a real function~$\f$.

It is important to note, however, that the inclusion $\FFR\subset\FRR$ is
proper, i.e., {\gr there exist functional relations in~$\FRR$ which do not
correspond, through the above identification, to any virtual function
of~$\FFR$}.

To prove that statement, let us consider the function $\Phi\:\RR\to\RR$ given
by:
$$
    \Phi(\x) = \cases{ 1 &, if $\x=2$ or $\x=-2$;\cr
                       0 &, otherwise.\cr}
$$
If $\Phi$ were a virtual function, then any representative sequence
$(\f_1,\f_2,\ldots)\in\SFR$ must end taking the value~$1$ at~$x=2$
and~$x=-2$, i.e., there would exist $n\in\N$ such that, for every $i>n$,
$\f_i(2)=\f_i(-2)=1$. But we would then have $\Phi(\pm2)=1$, and not
$\Phi(\pm2)=0$. (According to the notation established in Ref.~2, the virtual
$\pm2\in\RR$ is the class of the sequence $(-2,+2,-2,+2,\ldots)\in\SR$.)

It is also easy to verify that the identifications:
$$
    \FR = \KFR \subset \FFR \subset\FRR
$$
are compatible with the virtual extension of real functions {\gr considered
as relations between members of~$\R\subset\U$}, not as elements
of~$\FR\subset\U$.

\tit{III. Construction of Virtual Functions}

The aim of this section is to show that the identification $\FFR\subset\FRR$
is also compatible with the usual ``pointwise'' manipulations of real
functions: compositions and algebraic operations. This compatibility will
allow us to construct and represent virtual functions ``as if they were
real''.

We can consider the composition of real functions as a map which assigns to
each pair $\f,g\in\FR$ a third function $(\f\o\g)\in\FR$ given by
$(\g\o\f)(x)=\g[\f(x)]$. The domain of that composite function is the set
(eventually empty) of all values of~$x\in\dom\f$ for which $\f(x)\in\dom\g$.
Thus, the virtual extension of this map defines the {\gr composition of
virtual functions:\/} an operation that assigns to each pair of virtual
functions $\phi,\psi\in\FFR$ a third virtual function $(\psi\o\phi)\in\FFR$,
which will be called {\gr composite of~$\phi$ and~$\psi$}.

Since the composition of real functions is done ``pointwisely'', the
above definition of composition of virtual functions is perfectly
compatible with the identifications:
$$
    \FR = \KFR \subset \FFR \subset \FRR,
$$
and the domain of the composite $(\psi\o\phi)$ is the set (eventually empty)
of all values of~$\x\in\dom\phi$ for which $\phi(\x)\in\dom\psi$. For these
values we have $(\psi\o\phi)(\x)=\psi[\phi(\x)]$.

This shows that {\gr the subset $\FFR\subset\FRR$ of virtual functions is
closed with respect to the composition operation}. It is also not difficult
to verify that if a virtual function is {\gr inversible\/} then {\gr its
inverse is also a virtual function}.

Analogously, we can consider the addition of real functions as a map which
assigns to each pair $\f,g\in\FR$ a third function $(\f+\g)\in\FR$ given by
$(\f+\g)(x)=\f(x)+\g(x)$. The domain of this sum function is the intersection
(eventually empty) of the domains of $\f$ and~$\g$. Thus, the virtual
extension of this map defines an {\gr addition of virtual functions:\/} a map
which assigns to each pair of virtual functions $\phi,\psi\in\FFR$ a third
virtual function $(\phi+\psi)\in\FFR$, which will be called {\gr sum
of~$\phi$ and~$\psi$}.

Since the real functions addition is also done ``pointwisely'', the above
definition of virtual functions addition is also perfectly compatible with
the identifications:
$$
    \FR = \KFR \subset \FFR \subset \FRR,
$$
and the domain of the sum $(\phi+\psi)$ is the intersection (eventually
empty) of the domains of $\f$ and~$\g$. So we have
$(\phi+\psi)(\x)=\phi(\x)+\psi(x)$, for every $\x\in\dom(\phi+\psi)$.

This shows that {\gr the subset $\FFR\subset\FRR$ of virtual functions is
closed with respect to the addition operation}.

Proceeding this way, we can define the remainder algebraic operations
(subtraction, multiplication, division and exponentiation) with virtual
functions, verify its compatibility with the above identifications, and the
closure of the subset $\FFR\subset\FRR$ of virtual functions with respect to
these operations.

When we use Infinitesimal Calculus, most of the time we deal with functions
constructed by successive applications of compositions and algebraic
operations on a set of functions named {\gr elementary:\/} constant
functions, exponentials and logarithmics, direct and inverse trigonometric
functions.

The same procedure can be used to handle virtual functions: we take a set of
virtual functions considered ``elementary'' and, from them, we construct
others applying com\-posi\-tions and algebraic operations. That set of
``elementary virtual functions'' might include the exponentials,
logarithmics, trigonometrics (direct and inverse), considered virtual
functions according to the identification:
$$
    \FR=\KFR\subset\FFR,
$$
as well as the {\gr constant virtual functions}, i.e., those which satisfy
the virtual extension of the attribute (of real functions) ``to be
constant''. Those constant virtual functions are the classes in $\FFR$ which
can be represented by sequences $(\f_1,\f_2,\ldots)\in\SFR$ such that all
$\f_i$ are constant functions, but not necessarily equal among them. In other
words, a virtual function $\phi\:\RR\to\RR$ is constant when there exists a
virtual number $\a\in\RR$ such that $\phi(\x)=\a$, for every $\x\in\RR$. For
instance, $\k(\x)=\8$ defines a constant virtual function.

Having included these constant virtual functions among the elementary
functions, we can easily construct many non-real functions of~$\FRR$ which
will automatically belong to~$\FFR$. As an illustration, the expression:
$$
    \phi(\x) = {e^{(\x^2 - \8^2)}\over \cos(\pi \d \x)}
$$
specifies a virtual function $\phi\in\FFR$ which has been formally defined as
the class of the sequence $(\f_1,\f_2,\ldots)$ given by:
$$
    \f_n(x) = {e^{(x^2 - n^2)}\over \cos\left({\pi x\over n}\right)}.
$$
We can calculate $\phi(\x)$, for specific values of~$\x$, by direct
substitution:
$$
    \phi(0) = {e^{(0^2 - \8^2)}\over \cos(\pi \d 0)} = {e^{-\8^2}\over\cos0}
            = e^{-\8^2},
$$
or:
$$
    \phi(\8) = {e^{(\8^2 - \8^2)}\over \cos(\pi \d \8)} = {e^0 \over \cos\pi}
             = -1.
$$
However, this function~$\phi$ is not defined at some virtual numbers
$\x\in\RR$. For example, if $\x=\8/2$ then $\cos(\pi\d\x)=\cos(\pi/2)=0$, so
the denomination of the above fraction nullifies itself.

We will adopt, exactly as we do for real functions, the following convention:

{\gr Every time we define a virtual function by an expression, without
explicitly indicating its domain, one should understand that it is the set of
all virtual numbers for which that expression has mathematical meaning.}

This way, we will not need to explicitly state domains when this is not
relevant. For instance, let us consider the virtual function given simply by:
$$
    \psi(\x) = {\8 \over 1 + \8^2\x^2}.
$$
Since this expression is well defined for any $\x\in\RR$, it is understood
that:
$$
    \psi\:\RR\to\RR.
$$
Formally, $\psi\in\FFR$ is the class represented by the sequence
$(\g_1,\g_2,\ldots)$ of the real functions:
$$
    \g_n(x) = {n \over 1 + n^2 x^2}.
$$

In our ``set of elementary virtual functions'' we can also include other
members of~$\FFR$ defined by explicit presentation of a representative
sequence $(\f_1,\f_2,\ldots)$ of real functions. For example, we can take:
$$
    \f_n(x) = \cases{ n/2 &, if $|x| < 1/n$;\cr
                      0   &, if $|x| \ge 1/n$,\cr}
$$
and then make $\chi=\<\f_1,\f_2,\ldots>\in\FFR$. For this virtual function
$\chi\:\RR\to\RR$ we have:
$$
    \chi(\x) = \cases{ \8/2 &, if $|\x| < \d$;\cr
                        0   &, if $|\x| \EXT\ge \d$.\cr}
$$
This pair of clauses, nevertheless, is not sufficient to define $\chi(\x)$
for every $\x\in\RR$, since there are virtual numbers $\x$ which do not
satisfy either condition above.

\tit{IV. Continuity}

We will now extend the notion of continuity to virtual functions.

We can interpret continuity of a function $\f\in\FR$ at a point $x\in\R$ as a
relation between $\f$ and~$x$. So, the virtual extension of this relation
allows us to define {\gr continuity of virtual function $\phi\in\FFR$ at a
point $\x\in\dom\phi$}. If $\phi$ is the class of the sequence
$(\f_1,\f_2,\ldots)\in\SFR$, and $\x$ the class of $(x_1,x_2,\ldots)\in\SR$,
then $\phi$ is continuous at~$\x$ when there exists $n\in\N$ such that $\f_i$
is continuous at $x_i$, for every $i>n$.

Moreover, we will say simply that a virtual function {\gr is continuous\/}
when it is continuous at each point of its domain (as we do for real
functions). It is easy to see that this attribute of virtual functions is
equivalent to the virtual extension of the attribute ``is continuous'' for
real functions.

We know from Calculus that the composite of two continuous real functions is
also a continuous (real) function. Hence, we have from the Virtual Extension
Theorem (VET, Ref.~1) that {\gr the composite of two continuous virtual
functions is also a continuous virtual function}.

The VET also shows that {\gr a virtual function obtained from two other
continuous virtual functions by an algebraic operation is necessarily
continuous}.

The attribute of virtual functions ``to be constant'', defined in the
previous section, is the virtual extension of the attribute ``to be
constant'' applied to real functions. Since any constant real function is
continuous, the VET guarantees that {\gr every constant virtual function is
continuous, even when the constant is a virtual number}.

Still from the VET, we know that the exponential, logarithmic, and
trigonometric functions (direct and inverse), considered as virtual functions
by the identification $\FR=\KFR\subset\FFR$, are all continuous.

Therefore, all virtual functions constructed from these ``elementary'' ones,
through successive applications of compositions and algebraic operations, are
continuous.

For instance, the virtual functions $\phi$ and~$\psi$ viewed in the previous
section are continuous:
$$
    \phi(\x) = {e^{(\x^2 - \8^2)}\over \cos(\pi \d \x)}
    \qquad\hbox{and}\qquad
    \psi(\x) = {\8 \over 1 + \8^2\x^2}.
$$
This latter function $\psi:\RR\to\RR$ is continuous at each point of its
domain, including at real values of its argument, for which:
$$
    \psi(x) \= \cases{ \8 &, if $x=0$;\cr
                       0  &, if $x\ne0$.\cr}
$$
However, the virtual function $\chi\:\RR\to\RR$ (also seen in the previous
section) is not continuous at $\x=\d$, since the corresponding functions
$\f_n$ are not continuous at~$x=1/n$.

The continuity of a real function $\f$ at a point $x$ of its domain was
defined in the first part of this work$^2$ by the condition:
$$
    \a\=x \IMP \f(\a)\=\f(x).
$$
It is important to note that continuity of virtual functions, as defined in
this section, does not require that this condition hold. For example,
evaluating the function $\psi$ at the infinitesimal $\a=\sqrt{\d}$ we get:
$$
    \psi(\sqrt{\d}) = {\8 \over 1 + \8^2\d} = {\8 \over 1 + \8} \= 1,
$$
i.e., $\sqrt{\d}\=0$ but~$\psi(\sqrt{\d})\not\=\psi(0)$.

Analogously, defining {\gr uniform continuity\/} of a virtual function $\phi$
by the extension of the attribute ``to be uniformly continuous'', we will not
have necessarily that:
$$
    \a\=\b \IMP \phi(\a)\=\phi(\b).
$$

\tit{V. Derivation}

Our goal in this section is to define the derivative of virtual functions,
and show that the usual derivation algorithm of Calculus can be used ``as if
virtual functions were real''.

We can interpret the derivation process of functions in~$\FR$ as a map which
assigns a function $\f'\in\FR$ to each $\f\in\FR$, the domain of~$\f'$ being
the set (eventually empty) of points at which $\f$ is derivable. Thus, the
virtual extension of this map allows us to define the {\gr derivative\/}
$\phi'\in\FFR$ of any virtual function $\phi\in\FFR$. If $\phi$ is the class
of the sequence $(\f_1,\f_2,\ldots)\in\SFR$ then its derivative $\phi'\in\RR$
is the class of the sequence $(\f'_1,\f'_2,\ldots)\in\SFR$.

We will say that $\phi\in\FFR$ {\gr is derivable at~$\x\in\RR$} when
$\x\in\dom\phi'$. We can alternatively interpret the derivability of a
function $\f\in\FR$ at a point $x\in\R$ as a relation between $\f$ and~$x$,
and then define the {\gr derivability of a virtual function $\phi\in\FFR$ at
a point $\x\in\dom\phi$} through the virtual extension of this relation.
These two definitions are equivalent: if $\phi$ is the class of the sequence
$(\f_1,\f_2,\ldots)\in\SFR$, and $\x$ is the class
of~$(x_1,x_2,\ldots)\in\SR$, then $\phi$ is derivable at~$\x$ when there
exists $n\in\N$ such that $\f_i$ is derivable at~$x_i$ for every $i>n$. In
this case, $\phi'(\x)\in\RR$ is exactly the class of sequences which end
taking the values $\f'_i(x_i)$.

Furthermore, we will say that a virtual function {\gr is derivable\/} when it
is derivable at each point in its domain (as we do for real functions).
Clearly, this definition is also equivalent to the direct virtual extension
of the attribute ``to be derivable''.

In the first part of this work we saw a stronger condition of {\gr
differentiability\/} for real functions than simple derivability, and we
proved that this condition is equivalent to the demand for continuity of the
derivative at the considered point. We define the attribute ``to be
differentiable'' for virtual functions through the extension of the
corresponding attribute for real functions. So, the VET shows that {\gr a
virtual function is differentiable at a point $\x$ in its domain \IFF\ it is
derivable at~$\x$, and its derivative is continuous at this point}.

The VET also guarantees that the ``virtual derivation'' process applied to
real functions provides the same result of ``real derivation''.  Besides,
every constant real function is derivable and its derivative is the null
function. So every constant virtual function is derivable and its derivative
is the null function, even when the constant is a virtual number.  For
example, if $\k(\x)=\8$ for every $\x\in\RR$, then $\k'(\x)=0$ for every
$\x\in\RR$.

The {\gr derivation rules\/} are assertions about the derivation process
considered as a map from $\FR$ into~$\FR$. According to the preceding
definitions, the VET guarantees that all of them also hold for virtual
functions. For instance:
$$
    (\phi + \psi)' = \phi' + \psi'
$$
$$
    (\phi\psi)' = \phi'\psi + \phi\psi'
$$
$$
    (\phi\o\psi)' = (\phi'\o\psi)\psi'.
$$

Therefore, {\gr the usual derivation algorithm of Calculus can be used for
virtual functions exactly as if they were real}. Doing so, {\gr we should
treat virtual constants (present in an expression which defines a virtual
function) exactly as we do real constants (present in expressions which
define real function)}. Thus, if
$$
    \psi(\x) = \d \x^\8 + \8^2
$$
then
$$
    \psi'(\x) = \d \8 \x^{(\8-1)} + 0 = \x^{(\8-1)}.
$$

As an illustration involving many derivation rules, the chain rule inclusive,
we have:
$$
    \phi(\x) = {e^{(\x^2 - \8^2)}\over \cos(\pi\d\x)}
$$
implies:
$$\eqalign{
    \phi'(\x) & = { e^{(\x^2 - \8^2)}(2\x)\cos(\pi\d\x) +
                         e^{(\x^2 - \8^2)}\sin(\pi\d\x)\pi\d \over
                  \cos^2(\pi\d\x)}\cr
              & = { e^{(\x^2 - \8^2)} \over \cos^2(\pi\d\x)}
                  \left[ 2\x\cos(\pi\d\x) + \pi\d\sin(\pi\d\x)\right].\cr}
$$

In Ref.~2 we observed that the {\gr Leibnizian notation\/} for derivatives of
real functions can be interpreted in~$\RR$ as a quotient between
infinitesimals. This notation can also be generalized to represent the
virtual derivation process earlier defined: if a virtual function is
specified by an expression which provides the values of its {\gr dependent
virtual variable\/}~$\y$ from the values of its {\gr independent virtual
variable\/}~$\x$, then the symbol:
$$
    {d\y \over d\x}
$$
will be used to denote the dependent virtual variable of its derivative. For
instance, if
$$
    \y = {\8 \over 1 + \8^2\x^2}
$$
then:
$$
    {d\y \over d\x} = {-\8 \over (1 + \8^2\x^2)^2}(2\8^2\x)
                    = {-2\8^3\x \over (1 + \8^2\x^2)^2}.
$$

However, it is important to notice that the symbols $d\y$ and~$d\x$ do not
denote ``infinitesimal variations'' of $\y$ or~$\x$, since these variables
are virtual numbers themselves. So, if we are using the Leibnizian notation
for the derivative of a virtual function $\y=\phi(\x)$ then we should
interpret the ``quotient'' $(d\y/d\x)$ in a purely symbolic way, and not as a
quotient between infinitesimals (contrary to the $dx$ e~$dy$ present in the
Leibnizian notation for the derivatives of real functions$^2$). Therefore, we
write:
$$
    {d\y \over d\x} = \phi'(\x) \qquad\hbox{and not}\qquad
      {d\y \over d\x} \= \phi'(\x),
$$
whereas for the derivative of a generic real function $y=\f(x)$ we write:
$$
    {dy \over dx} \= \f'(x) \qquad\hbox{and not}\qquad
      {dy \over dx} = \f'(x).
$$

The derivative of a real function $\f$ at a point $x$ in its domain was
defined in the first part of this work as being the real number $\f'(x)$
which satisfies the condition:
$$
    \a\@ x \IMP \f'(x) \= { \f(\a) - \f(x) \over \a - x }.
$$
Nevertheless, the derivative of a virtual function at a point in its domain,
as defined in this very section, does not necessarily satisfy this condition.
For example, calculating the derivative of the above function $\y=\phi(\x)$
at~$\x=0$ we get:
$$
    \phi'(0) = \left.{d\y \over d\x}\right|_{\x=0}
        = {-2\8^3 0 \over (1 + \8^2 0^2)^2} = 0,
$$
(which was expected, since $\phi$ is an even function) whereas, for $x=0$
and~$\a=\d$, we have:
$$
\eqalign{ { \phi(\a) - \phi(x) \over \a - x }
    &= {{\8 \over 1 + \8^2\d^2} - {\8 \over 1 + \8^2 0^2}\over\d - 0 }\cr
    &= { {\8 \over 1 + 1} - {\8 \over 1 + 0} \over\d}\cr
    &= \8 \left( {\8 \over 2} - \8 \right)\cr
    &= -{\8^2 \over 2}.\cr}
$$
So $\d\@0$ but
$$
    \phi'(0) \not\= { \phi(\d) - \phi(0) \over \d - 0 }.
$$

\tit{VI. Integration}

We will now define the integral of virtual functions. The integration theory
on the real line due to Riemann is sufficient for our purposes. So, from here
on, integrability and integrals of real functions between two real numbers
should be understood in this sense.

For any real number $a$, we will consider that a real function is
Riemann-integrable ``between $a$ and~$a$'' when $\f$ is defined at~$a$, and
in this case we shall have:
$$
    \int_a^a\f(x)\,dx = 0.
$$
Besides, we will not require that the lower integration limit be less than
the upper one. We will always suppose that:
$$
    \int_a^b\f(x)\,dx = - \int_b^a\f(x)\,dx.
$$

We can consider the integrability of a real function $\f\in\FR$ between two
points $a,b\in\R$ as a relation between $\f$ and those two real numbers, and
so define the {\gr integrability of a virtual function $\phi\in\FFR$ between
two virtual points\/} $\a,\b\in\RR$ through the virtual extension of this
relation. In other words, if $\phi$ is the class of the sequence
$(\f_1,\f_2,\ldots)\in\SFR$, $\a$ is the class of~$(a_1,a_2,\ldots)\in\SR$,
and $\b$ is the class of~$(b_1,b_2,\ldots)\in\SR$, then $\phi$ is integrable
between ~$\a$ and~$\b$ when there exists $n\in\N$ such that, for every $i>n$,
the function $\f_i$ is integrable between ~$a_i$ and~$b_i$.

Moreover, we will merely say that a virtual function {\gr is integrable\/}
when it is integrable between any pair of virtual numbers in its domain.

If $\phi\in\FFR$ is integrable between $\a$ and~$\b$, then we define the {\gr
virtual integral}
$$
    \int_\a^\b\phi(\x)\,d\x
$$
through the virtual extension of the map which assigns the number
$$
    \int_a^b\f(x)\,dx
$$
to each triple $(\f,a,b)$ with $\f$ integrable between $a$ and~$b$. The use
of a Greek letter as {\gr integration variable\/} indicates the integral is
virtual.

It is easy to see that if $\phi=\<\f_1,\f_2,\ldots>$ is integrable between
$\a=\<a_1,a_2,\ldots>$ and $\b=\<b_1,b_2,\ldots>$ then the integral of $\phi$
between ~$\a$ and~$\b$ is the class of sequences which end taking the values:
$$
    \int_{a_i}^{b_i}\f_i(x)\,dx.
$$
Thus, it should be clear that if $\f$ is a real function integrable between
$a$ and~$b$ then:
$$
    \int_a^b\f(\x)\,d\x = \int_a^b\f(x)\,dx,
$$
i.e., {\gr it does not matter whether we use a real (Latin letter) or a
virtual (Greek letter) integration variable when the integrand and the
integration limits are real.}

We will say that a virtual function {\gr is defined between two virtual
numbers\/} when those three objects satisfy the virtual extension of the
corresponding relation between a real function and two real numbers.
Analogously, we will say that a virtual function {\gr is con\-tinu\-ous
between two virtual numbers\/} when those three objects satisfy the virtual
extension of the corresponding relation between a real function and two real
numbers.

We know that if a real function is defined and continuous between two real
numbers then it is integrable between those two numbers. So, the VET
guarantees that {\gr every virtual function defined and continuous between
two virtual numbers is integrable between those two numbers}.

Therefore, the symbols:
$$
    \int_0^\d {\8 \over 1 + \8^2\x^2}\,d\x
$$
and
$$
    \int_{-\8^2}^{\8^\x} e^{\d\t}\,d\t
$$
are perfectly defined (they are just two particular examples). The last one
should be understood as a virtual variable which depends on~$\x$ through a
virtual function implicitly indicated. That function is the class of
sequences which end equals to $(\f_1,\f_2,\ldots)\in\SFR$ given by:
$$
    \f_n(x) = \int_{-n^2}^{n^x} e^{nt}\,dt.
$$

It is also known that the sum and the product of two integrable real
functions are also integrable. So, the VET guarantees that {\gr the sum and
the product of two integrable virtual functions are also integrable}.

In the same way, we can use the VET to show that virtual integrals satisfy
many of the properties expected from an integral, as, for example, the {\gr
additivity with respect to the integration interval:}
$$
    \int_\a^\ga\phi(\x)\,d\x = \int_\a^\b\phi(\x)\,d\x
          + \int_\b^\ga\phi(\x)\,d\x.
$$
We will not list all those properties here, since it is more convenient to
use the VET (which is always at hand) whenever these properties are required.
However, the basic fact from Differential and Integral Calculus which relates
derivatives and integrals deserves a more detailed discussion.

\tit{VII. The Fundamental Theorem of Calculus}

We will say that $\psi\in\FFR$ is a {\gr primitive\/} of $\phi\in\FFR$ when
$\psi'=\phi$. The set of all primitives of a virtual function will be called
its {\gr indefinite integral\/} (as we do for real functions).

The VET guarantees that {\gr two primitives of a continuous virtual function
$\phi\:\RR\to\RR$ differ by a (virtual) additive constant}. Thus, we will
represent the indefinite integral of that function by:
$$
    \int \phi(\x)\,d\x = \psi(\x) + \k,
$$
where $\psi\:\RR\to\RR$ is a particular primitive of $\phi$ and $\k\in\RR$ a
{\gr generic constant}, being understood that we are referring to the set of
virtual functions which can be written in that manner.

According to those definitions and conventions, the VET shows that the
integration techniques of traditional Calculus extend to virtual functions,
as well as the traditional notations associated to them, for which it is
enough to change from Latin to Greek letters.

For instance, to calculate
$$
    \int {\8 \over 1 + \8^2\x^2}\,d\x
$$
we make the substitution of variables $\m=\8\x$, so $d\m=\8d\x$ and:
$$
    \int {\8 \over 1 + \8^2\x^2}\,d\x = \int {d\m \over 1 + \m^2}
    = \arctan\m + \k = \arctan(\8\x) + \k.
$$

Applying the VET to the first form of the Fundamental Theorem of Calculus we
get:

{\gr If $\phi\in\FFR$ is defined and continuous between $\a$ and~$\b$, then,
for every $\x$ between $\a$ and~$\b$:}
$$
    {d\over d\x} \int_\a^\x\phi(\t)\,d\t = \phi(\x).
$$

In other words, the virtual function defined by:
$$
    \psi(\x) = \int_\a^\x\f(\t)\,d\t
$$
is a primitive of~$\phi$. In particular, {\gr every continuous virtual
function $\phi\:\RR\to\RR$ admits a primitive}.

The VET also provides the ``virtual second form'' of the Fundamental Theorem
of Calculus:

{\gr If $\phi\in\FFR$ is defined and continuous between $\a$ and~$\b$, and
$\psi$ is a primitive of~$\phi$, then:}
$$
    \int_\a^\b\phi(\x)\,d\x = \psi(\b) - \psi(\a).
$$

Some examples:
$$
    {d\over d\x} \int_{-\8^2}^{\8^\x} e^{\d\t}\,d\t
    = e^{\d\8^\x} {d\over d\x} \left(\8^\x\right)
    = e^{\8^{\x-1}} \8^\x \ln\8.
$$
$$
    \int_0^\d {\8 \over 1 + \8^2\x^2}\,d\x = \left[\arctan(\8\x)\right]_0^\d
    = \arctan(\8\d) - \arctan(\8 0) = \arctan 1 - \arctan 0 = {\pi\over 4}.
$$

In Sec.~VI we saw that it does not matter whether we use a real (Latin
letter) or a virtual (Greek letter) integration variable when the integrand
and the integration limits are real. When this is not the case, the symbol
$$
    \int_\a^\b\phi(x)\,dx
$$
has not yet been defined.

We will say that a virtual integral
$$
    \int_\a^\b \phi(\x)\,d\x
$$
is {\gr reducible\/} when it is near some real number. In this case, we will
use the symbol
$$
    \int_\a^\b\phi(x)\,dx
$$
to denote this real number, which will be called {\gr reduced integral\/} of
function $\phi$ between $\a$ and~$\b$. This symbol will not be used if the
virtual integral is not reducible, so an integral with a real (Latin letter)
integration variable is always a real number, if it exists.

For any reducible integral, we have:
$$
    \int_\a^\b\phi(x)\,dx \= \int_a^b\phi(\x)\,d\x,
$$
but, in general, it is not true that:
$$
    \int_\a^\b\phi(x)\,dx = \int_a^b\phi(\x)\,d\x.
$$

For instance:
$$
    \int_1^\8 {d\x\over\x^2} = \left[{-1\over\x}\right]_1^\8 = -\d + 1,
$$
so
$$
    \int_1^\8 {dx\over x^2} = 1,
$$
and
$$
    \int_1^\8 {dx\over x^2} \not= \int_1^\8 {d\x\over \x^2}.
$$

Analogously:
$$
    \int_\d^1 {d\x\over\sqrt\x} = \left[2\sqrt\x\right]_\d^1 = 2 - 2\sqrt\d,
$$
so
$$
    \int_\d^1 {dx\over\sqrt x} = 2,
$$
hence:
$$
    \int_\d^1 {dx\over\sqrt x} \not= \int_\d^1 {d\x\over\sqrt\x}.
$$

\tit{VIII. Virtual Sequences}

The virtual functions whose domain is the set $\NN\subset\RR$ will be called
{\gr virtual sequences}, and they will be generically denoted by~$(\a_\n)$,
where the {\gr virtual index\/} $\n$ ranges over the set~$\NN$ of all virtual
natural numbers. Since we can calculate
$$
    \sum_{n=1}^k a_n
$$
for every $k\in\N$, the sum
$$
    \sum_{\n=1}^\k \a_\n
$$
is always a well defined virtual number, for any virtual sequence $(\a_\n)$
and any $\k\in\NN$. Thus, we can add infinitely many virtual numbers:
$$
    \sum_{\n=1}^\8 \a_\n
$$
without worrying about ``technical'' convergence questions.

As an example of application, we can use Riemann sums of infinitely many
infinitesimal terms to integrate a real function $\f$ defined between two
real numbers $a$ and~$b$. To do so, it is enough to consider an {\gr
infinitely fine partition\/} of that interval as a virtual sequence $(\a_\n)$
such that both $\a_1=a$ and $\a_\8=b$, as well as:
$$
    \hbox{$\a_{\n+1}\=\a_\n$, \qquad for every $\n\in\NN$.}
$$
As an illustration, the virtual sequence:
$$
    \a_\n = {\n\over\8}
$$
is an infinitely fine partition of the real interval $[0,1]$ according to
that definition.

At this point, the general mechanism which allows us to transport the
constructions and techniques of Calculus to virtual functions and sequences
should be clear. After having understood this mechanism, it is easy to see
that we could have defined {\gr virtual subsets\/} of the real line, as well
as the notion of {\gr virtual relation}, from the very beginning. Through
obvious identifications, these objects could have been considered,
respectively, as subsets of~$\RR$ and relations between virtual numbers.

We would then have, for instance, that the set of all virtuals between $\a$
and~$\b$ would be a virtual subset of~$\RR$, for any pair $\a,\b\in\RR$. More
generally, the domain of any virtual function would always be a virtual set.

That general mechanism also allows us to construct {\gr virtual differential
equations}, which can have virtual functions as coefficients or initial
conditions. The {\gr unknown\/} of such an equation is a virtual function,
and we can talk about general or particular {\gr virtual solutions}.

It is not difficult to see that we can also reformulate the Multivariable
Calculus as we did in the first part of this work for just one variable, and
so ``virtualize'' many other notions, like differential partial equations and
their respective solutions, for example.

In the context of Functional Analysis, or of Linear Algebra in general, the
idea of virtual sequence of vectors can be used to define {\gr infinite
linear combinations}, which leads to the notion of {\gr virtual base\/} of a
vector space. Perhaps this concept reduces the ``tech\-nic\-alities''
involved in the representation of functions by functional series, like
Fourier's, for instance.

We do not intend to explore all those possibilities here, but our hope is to
have clearly indicated a way which might enlarge Infinitesimal Calculus
usefulness and reach.

\tit{References}

\item{$^1$} S.~F.~Cortizo: ``Virtual Extensions'', to appear (1995).

\item{$^2$} S.~F.~Cortizo: ``Virtual Calculus---Part I'', to appear (1995).

\item{$^3$} P.~A.~M.~Dirac: {\it The Principles of Quantum Mechanics\/}
(Oxford at the Clarendon Press, 1947), 3rd ed., pp.~58--61.
\bye